\documentclass[reqno, 11pt]{amsart}
\usepackage{fullpage}
\usepackage{stmaryrd}
\usepackage{amssymb}
\usepackage{mathrsfs}
\usepackage{amsfonts,amssymb}
\usepackage{cite}

\newtheorem{thm}{Theorem}[section]

\newtheorem{prop}[thm]{Proposition}

\newtheorem*{conjecture*}{Conjecture}
\newtheorem*{thm*}{Theorem}

\theoremstyle{definition}

\theoremstyle{remark}

\numberwithin{equation}{section}

\newcommand{\Aut}{{\operatorname{Aut}}}

\newcommand{\Irr}{{\operatorname{Irr}}}

\newcommand{\Syl}{{\operatorname{Syl}}}

\newcommand{\Sp}{{\operatorname{Sp}}}

\newcommand{\GL}{{\operatorname{GL}}}

\begin{document}

\title{}

\author{}

\author{}

\address{}

\address{}

\makeatletter

\makeatother

%\keywords{Brauer character, Sylow $p$-subgroup, derived subgroup}
\date{}

%\dedicatory{}

% at present the "communicated by" line appears only in ERA, PROC and JAG
%\commby{}

\title{Bounding the number of $p'$-degrees from below}

\author[T. M. Keller]{Thomas Michael Keller}
\address{Department of Mathematics, Texas State University, San Marcos}
\email{keller@txstate.edu}

\author[Y. Yang]{Yong Yang}
\address{Department of Mathematics, Texas State University, San Marcos}
\email{yang@txstate.edu}


%\author{Yong Yang}

%\address{Key Laboratory of Group and Graph Theories and Applications, Chongqing University of Arts and Sciences, Chongqing, China and Department of Mathematics, Texas State University, 601 University Drive, San Marcos, TX 78666, USA.}

\begin{abstract}
Let $G$ be a finite group of order divisible by a prime $p$
and let $P\in\Syl_p(G)$.
We prove a recent conjecture by Hung stating that
$|\Irr_{p'}(G)|\geq \frac{\exp(P/P')-1}{p-1}+2\sqrt{p-1}-1.$
Let $a\geq 2$ be an integer and suppose that $p^a$
does not exceed the exponent of the center of $P$.
We then also show that the number of conjugacy classes of elements of $G$ for which $p^a$ is the exact $p$-part of their order is at least $p^{a-1}$.

%We obtain a result regarding the existence of small orbits for self-dual modules, thus generalizing a theorem of T. Berger for symplectic modules. We also strengthen a result of Marchi.
\end{abstract}

\maketitle

\Large

\section{Introduction}

Let $G$ be a finite group and $p$ a prime which divides the order of $G$. The number $k(G)$ of conjugacy classes of $G$ is a fundamental invariant in group and representation theory, and bounding $k(G)$ from below  in terms of $p$ is a question that has been studied intensely in the past two decades.
One reason that it has attracted so much interest is that it is related to a well-known problem of Brauer
\cite[Problem 21]{Bra63}. A (still open) conjecture of H{\'e}thelyi and K{\"u}lshammer \cite{HK00} is that for any $p$-block $B$ of any finite group $G$, the number of $k(B)$ of complex irreducible characters in $B$ is $1$ or is at least $2\sqrt{p-1}$. In particular, this implies $k(G) \geq 2\sqrt{p-1}$ for all $G$ and $p$.

Proving $k(G)\geq 2\sqrt{p-1}$ for all $G$ and $p$ turned out to be a difficult problem. Building on a series of relevant works by H{\'e}thelyi-K{\"u}lshammer{\cite{HK00,HK03}, Malle~\cite{Mal06}, Keller~\cite{Kel09}, and H{\'e}thelyi-Horv{\'a}th-Keller-Mar{\'o}ti~\cite{HHKM11}, the conjecture was confirmed in ~\cite{Mar16}. Further work~\cite{MS20} showed that there is some constant $c > 0$ such that when the square of the prime divides the order of the group, then the bound $k(G) \geq cp$ holds. The problem has been taken into various other directions, mostly replacing $k(G)$ by smaller quantities or improving the lower bound by imposing additional hypotheses.
In this paper we improve the lower bound under the additional hypothesis
that the exponent of the Sylow $p$-subgroups is greater than $p$.

Let $a$ be a non-negative integer. Throughout this paper we use $k_G(p^a)$ to denote the number of conjugacy classes of those elements $g$ of $G$ satisfying $|g|_p=p^a$, where $|g|_p$ is the $p$-part of the order of $g$.
Another way to say this is that if $\Omega$ is the set of elements $g$ of $G$
with $|g|_p=p^a$, then $k_G(p^a)$ is the number of orbits of $g$ acting via conjugation on $\Omega$.

With this we prove the following.

%We define the \textit{$p$-regularity level} of an element $g\in G$ to be $\log_p(|g|_p)$, where $|g|_p$ is the $p$-part of the order of $g$. Throughout this paper we use $k_{pr-a}(G)$ to denote the number of conjugacy classes of $G$ with $p$-regularity level $a$.

%We bound the number of $p$-regularity level $a$ conjugacy classes of all finite groups in terms of $a$ and $p$.

\begin{thm}
  \label{thm:pra1}
  Let $G$ be a finite group of order divisible by $p$, where $p$ is a prime. Let $P\in\Syl_p(G)$ and $a$ an integer such that $a\leq \log_p(\exp(\mathbf{Z}(P)))$, Then
  \[k_{G}(p^a)\geq p^{a-1}.\]
\end{thm}

We note that the bound in Theorem \ref{thm:pra1} is sharp,
as the following example shows. Let $P$ be cyclic of order $p^b$, and let $A\leq \Aut(G)$ such that $|A|=p-1$. Then
$G=AP$ is a Frobenius group, and for any $1\leq a\leq b$ we see that $k_{G}(p^a) = (p^a-p^{a-1})/(p-1)= p^{a-1}$.\\

\indent As a consequence of our work to prove Theorem \ref{thm:pra1} we can also settle a recent conjecture by N.\,N. Hung presented in \cite[Conjecture 1.4]{hung2023} which strengthens the aforementioned lower bounds on $k(G)$ considerably. In
\cite[Theorem 1.5]{hung2023} he proved the conjecture for $p=2$ using $p$-rationality arguments and noted in his Remark 1.7 that
by using the McKay conjecture for $p=2$ (which was the only proved case when he wrote his paper) the problem ``is reduced to showing that the
conjugacy class number $k(AM)$ of the semidirect product of an odd-order group
$M$ acting on an abelian 2-group $A$ is at least exp($A$). However, even in this much
simpler situation, we are not aware of any proof that does not use the idea of
2-rationality level (...)". Here we will do just that: find an elementary proof of this
using orbit counting arguments, and not just for $p=2$, but for all primes $p$ thanks to the recent completion of the proof of the McKay conjecture for all primes $p$. We thus prove:

\begin{thm}
  \label{thm:hung}
  Let $G$ be a finite group of order divisible by $p$, where $p$ is a prime. Let $P$ be a Sylow $p$-subgroup of $G$. Furthermore, let $a$ and $b$ be positive integers such that $a-b$ is minimal and $p-1=ab$. Then
  \[|\Irr_{p'}(G)|\geq \frac{\exp(P/P')-1}{p-1}+2\sqrt{p-1}-1.\ \ \ \ (1)\]
  Moreover, there exists a universal constant $C$ such that if $p>C$, then we even have
  \[|\Irr_{p'}(G)|\geq \frac{\exp(P/P')-1}{p-1}+a+b-1.\ \ \ \ (2)\]
\end{thm}

Observe that the first statement in Theorem \ref{thm:hung} proves Hung's \cite[Conjecture 1.4]{hung2023}, and the second statement is even stronger. It is likely that the second statement is true for all primes $p$, but currently we can only
prove it for large $p$. \\
\indent In view of \cite[Conjecture A]{cinarci-keller2023} we conjecture that the second statement of Theorem \ref{thm:hung} is true for all primes $p$ (i.e. one can choose $C=1$).
Also, it might be interesting to determine when equality holds. When $\exp(P/P')=p$, then the first and last term in (1) and (2) cancel out and equality has been studied, see e.g.
\cite{cinarci-keller2023} or \cite{malle-maroti}. So it might
be interesting here to examine when equality in (1) or (2) can occur under the assumption that $\exp(P/P')\geq p^2$. For cyclic 2-groups both (1) and (2) are sharp, but this is not the
most exciting example. For $p=3$, we can look at the example following the statement of Theorem \ref{thm:pra1}, and this shows that (2) is sharp, and while the right hand side of (1) is not an integer value in this case, since we know that the
left-hand side of (1) is an integer, we can view (1) sharp as
well. These seem to be the only cases where (1) and (2) can provide a sharp lower bound when the exponent of the Sylow $p$-subgroup exceeds $p$.\\

\section{CONJUGACY CLASSES OF ELEMENTS WITH VARYING $p$-PARTS}

\begin{prop}\label{prop4}
 Let $p$ be a prime number. Then for every permutation group $G\leq S_n$ of degree $n>1$ with $p$ not dividing the order of $G$, we have $|G|<p^{n-1}$.
\end{prop}
\begin{proof}
This is \cite[Proposition 4]{PL}.
\end{proof}

\begin{thm}\label{theorem 5}
Let $G\leq \GL(d,p)$ be a primitive linear group such that $p\nmid |G|$. Then $|G|\leq p^d\cdot d\cdot \log_2p$, unless $p=7,\ d=4$ and $G\cong \Sp(4,3)$ or $\Sp(4,3)\times Z_3$, or $p=3,\ d=4$, $|G| = 4 \cdot 5 \cdot 2^5$, and $G$ has exactly two orbits on $V$.
\end{thm}
\begin{proof}
This is the main result of ~\cite{AT}, see also \cite[Theorem 2.3]{Kel09}.
\end{proof}

\begin{prop}
\label{coprimeorder}
    Let $A$ be a subgroup of $\GL(k, p)$ that has order coprime to $p$, then $|G| \le p^{k-1} \cdot (p^k-1)$.
\end{prop}
\begin{proof}
By Maschke's theorem, the action of $A$ on $V$ is completely reducible. If it is not irreducible, say $V=V_1\times V_2$ is a nontrivial factorization where $|V_1|=p^{k_1}$ and $|V_2|=p^{k_2}$, then induction yields \[|A|\leq|A/C_A(V_1)|\cdot |A/C_A(V_2)|\leq p^{k_1-1} \cdot (p^{k_1}-1) \cdot p^{k_2-1} \cdot (p^{k_2}-1) \le p^{k-1} \cdot (p^k-1).\] Hence we may assume that $A$ acts irreducibly on $V$, that is, $A\leq\GL(k,p)$ is an irreducible linear group.

If $A$ is an imprimitive linear group, then we can embed $A$ into the wreath product of a primitive linear group $B\leq \GL(d,q)$ and a permutation group $H\leq S_n$, where $p\nmid|B|,\ p\nmid|H|$ and $k=dn,\ n\geq2$. If $A$ is a primitive linear group, then let $B=A,\ n=1$. Combining Proposition \ref{prop4} and Theorem \ref{theorem 5}, we obtain
$$|A|\leq|B|^n|H|\leq(p^d\cdot d\cdot\log_2p)^np^{n-1},$$
with the noted exceptions $p=7,\ d=4$ or $p=3,\ d=4$.

So our desired inequality $|A|\leq(p^{dn}-1)(p^{dn-1})$ holds whenever $d\cdot \log_2p\leq p^{d-1}-1$, which is true except for $d=1$ and for $d=2,\ p=2,3,$ or $5$. However, for $d=1$ we obviously have $|B|\leq p-1$, hence $|A|\leq (p-1)^n p^{n-1}$. For $d=2,\ p=2$, we have $|B|\leq 3$, $|A|\leq 3^n \cdot 2^{n-1}<2^{2n-1} (2^{2n}-1)$; for $d=2,\ p=3$ we obtain $|B|\leq16,\ |A|\leq 16^n \cdot 3^{n-1}<3^{2n-1} (3^{2n}-1)$; and for $d=2,\ p=5$ we have $|B|\leq 96$, $|A|\leq 96^n \cdot 5^{n-1}<5^{2n-1} (5^{2n}-1)$. For $p=7,\ d=4,\ B \leq  \Sp(4,3) \times Z_3$, we obtain $|A|\leq155520^n\cdot7^{n-1}<7^{4n-1}(7^{4n}-1)$. Finally, for $p=3,\ d=4$, $|B| = 4 \cdot 5 \cdot 2^5$, and thus $|A|\leq 640^n\cdot3^{n-1}<3^{4n-1}(3^{4n}-1)$.
\end{proof}

%Consider the action of the semidirect product of $H$ with $\bZ(P)$ where the order of $H$ is coprime with $p$,  then the number of classes at level $a$ where $a \geq 2$ is
%\begin{equation}\label{1}
    %\frac{[(p^a)^m-(p^{a-1})^m]\cdot p^{2l}}{p^{2n-2}\cdot (p-1)}=\frac{p^{(a-3)m}[p^m-1]\cdot p^{2m}\cdot p^{2l}}{p^{2n-2}\cdot (p-1)}=\frac{p^{(a-3)m}[p^m-1]\cdot p^2}{ (p-1)}\ge p^{a-1}
%\tag{*}\end{equation}($m$ is the number of terms of level $\ge$ $a$, where $l$ is the number of terms of level $<$ $a$, clearly $n = m + l$).

%Using counting argument we may prove that the number of orbits of $H$ on the $a$-level conjugacy classes is $\ge$ $ p^{a-1}$ if $k_l\ge 2$.

\begin{thm}
  \label{orbits}

Let $P>1$ be an abelian $p$-group and let $G$ be $p'$-group which acts faithfully (via automorphisms) on $P$. Let $a \in \mathbb{Z} $ such that $1 \le a \le \log_p(\exp{P})$ and let $\Gamma=\Gamma_a$ be the set of elements of order $p^a$ in $P$. Write $n(G, \Gamma)$ for the number of orbits of $G$ on $\Gamma$. Then
$n(G, \Gamma)\ge p^{a-1}$.
\end{thm}

\begin{proof}
As is common, for a non-negative integer we write
$\Omega_i(P)=\langle x\in P\ |\ x^{p^i}=1\rangle$. \\
The case that $a=1$ is trivial, so let $a\geq 2$ be an integer. We show that for this $a$ a counterexample cannot exist. So let $G$ and $P$ be as in the
theorem such that  $a\leq \log_p(\exp(P))$ and such that $n(G,\Gamma) < p^{a-1}$. We may choose $G$ and $P$
with these properties such that $|GP|$ is minimal.\\
Since $\Omega_i(P)$ is normal in $GP$ for all $i$ and by the minimality of $GP$, we may assume that $P=\Omega_a(P)$.\\
Next observe that by \cite[page 175]{gorenstein} we know that there exists a positive integer $n$ and $G$-invariant subgroups $P_i$ of $P$ ($i=1,\dots, n$) such that
\[P=P_1\times\dots\times P_n,\]
and no $P_i$ can be decomposed further as the
direct product of nontrivial $G$-invariant subgroups.
The $P_i$'s are then called indecomposable. Clearly we may assume that the exponent of $P_1$ is $p^a$.
If $G$ acts trivially on $P_1$, then clearly
$n(G,\Gamma)\geq p^a$ and we are done. So
we may assume that $G$ does not act trivially on
$P_1$.
Then by our minimality assumption we may assume that
$P=P_1$. Since $P$ is indecomposable, from the well-known Zassenhaus decomposition (also known as Fitting's lemma, see \cite[Theorem 4.34]{isaacsfinitegroups}) we know that $P=[P, G]$, and by minimality $G$ acts faithfully on $P$. We now invoke
\cite[Theorem 5.2.2]{gorenstein} to conclude that $P$ is homocyclic, i.e., a direct product of $m$ copies of cyclic groups of order $p^a$. Thus $G$ acts faithfully on $P/\Phi(P)$ and $|G| \leq p^{m-1}(p^m-1)$ by Proposition ~\ref{coprimeorder}.

Then the number of orbits of $G$ on $\Gamma$ is at least
%Now let $k$ be the minimal number of generators of the elementary abelian group $V=P/ \Omega_{a-1}(P)$, and choose
%$x_i\in P$ ($i=1,\dots ,k$) such that $V=\langle \bar{x_1}, \dots , \bar{x_k}\rangle$, where $\bar{x_i}=x_i\Omega_{a-1}(P)$.
%Then $G$ acts on $V$, and if we put $Q=\langle x_1, \dots , x_k\rangle$. Furthermore, all $x_i$ have order $p^a$, and $Q$ is homocyclic
%$of type $(p^a)^k$. \\
%Now by coprime action (Zassenhaus decomposition) there exists a $G$-invariant $R\leq P$ such that $P=Q\times R$. By the
%minimal choice of $GP$ we may assume that $R=1$.

%Consider the action of the semidirect product of $H$ with $\bZ(P)$ where the order of $H$ is coprime with $p$,

\begin{equation}\label{1}
    \frac{(p^a)^m-(p^{a-1})^m}{p^{m-1}\cdot (p^m-1)}=\frac{p^{(a-1)m}(p^m-1)}{p^{m-1}\cdot (p^m-1)}=p^{(a-2)m+1}\ge p^{a-1}
\tag{*}\end{equation}.

But then by $(*)$  we see that $GP$ is not a counterexample. This
contradiction completes the proof.
\end{proof}

%\begin{thm}
  %\label{thm:pra}
  %Let $G$ be a finite group of order divisible by $p$, for $p$ prime. Then for $2\leq a\leq \log_p(\exp(\mathbf{Z}(P)))$ where $P\in\Syl_p(G)$ we have,   \[k_{pr-a}(G)\geq p^{a-1}.\]
%\end{thm}
We now prove Theorem \ref{thm:pra1}.
%\gap
\begin{proof}
We will find suitable conjugacy classes of $G$ from elements in $\mathbf{Z}(P)$. Suppose that $x$ and $y$ are in $\mathbf{Z}(P)$  and they are conjugate in $G$, then they are also conjugate in $N_G(P)$ (cf. \cite[Lemma 5.12]{isaacsfinitegroups}). Thus we may assume that $P$ is normal in $G$. By the Schur-Zassenhaus theorem, $P$ has a complement $H$ in $G$. Since $\mathbf{Z}(P)$ is characteristic in $P$, $\mathbf{Z}(P)$  is normal in $G$. We now consider the group action of $H/C_H(\mathbf{Z}(P))$ on $\mathbf{Z}(P)$ and the result follows by  Theorem ~\ref{orbits}.
\end{proof}

%We first note that for $P\in \Syl_p(G)$ we have the following factorization $\bZ(P)\cong \mathbb{Z}_{p^{k_1}}\times \dots \times \mathbb{Z}_{p^{k_n}}$.

\section{HUNG'S CONJECTURE}
In this section we prove Theorem \ref{thm:hung}.
\begin{proof}
Let $G$ be a counterexample to the theorem of minimal order.
Since the McKay conjecture now is a theorem as recently announced by B. Sp\"{a}th, we can use it and thus may assume that $P$ is normal in $G$.
By minimality we then may further assume that $P'=1$, i.e., $P$ is abelian, and then it is well-known that $\Irr_{p'}(G) = \Irr(G)$ and hence $|\Irr_{p'}(G)| = k(G)$.
Therefore we have to show that
\[k(G)\geq \frac{\exp(P)-1}{p-1}+2\sqrt{p-1}-1,\]
and
\[k(G)\geq \frac{\exp(P)-1}{p-1}+a+b-1\]
for large $p$.
Also, by Schur-Zassenhaus theorem $P$ has a complement $H$ in $G$, and by minimality we may assume that $H$ acts faithfully on $P$. So we have $G=HP$ where $H$ acts coprimely and
faithfully on the abelian group $P$.\\
Now let $2\leq a \leq \log_p(\exp(P))$. Then by Theorem \ref{orbits}
we have $n(G, \Gamma_a)\ge p^{a-1}$, where
$\Gamma_a$ denotes the set of elements of order $p^a$ in $P$, and
$n(G, \Gamma)$ is the number of orbits of $G$ on $\Gamma_a$.
Since orbits of $H$ on $P$ are conjugacy classes of $HP=G$, from
orbits of $H$ on elements of order $\geq p^2$ on $V$ we thus get
at least
\[\sum\limits_{a=2}^{\log_p(\exp(P))} p^{a-1}=p+p^2+\cdots+ p^{\log_p(\exp(P))-1} =  \frac{\exp(P)-1}{p-1}-1  \ \ \ (**)\]
conjugacy classes of $G$.\\
Next observe that $Q:=\Omega_1(P)$ (which by definition is the set
of elements of order $p$ or 1 of $P$) is characteristic in $G$, and
$H$ acts faithfully on $Q$ by \cite[Theorem 5.2.4]{gorenstein}. Furthermore, no element in $HQ$ has
order divisible by $p^2$, which implies that the the set of conjugacy
classes of $HQ$ is disjoint from the set of conjugacy classes counted in $(**)$, Thus
\[k(G)\geq \frac{\exp(P)-1}{p-1}-1 + k(HQ).\]
Now by the main result of \cite{Mar16} we conclude that $k(HQ)\geq 2\sqrt{p-1}$ which gives
us the first statement of the theorem, and \cite[Theorem C]{cinarci-keller2023} implies the second
statement of the theorem, completing the proof.
\end{proof}

\section*{Acknowledgements}
 Yang was partially supported by a grant from the Simons Foundation (\#918096, YY).

\end{document}